\documentclass{amsart}

\title{Bornological modifications of hyperspace topologies}
\author{Tom Vroegrijk}
\date{}

\newcommand{\cl}{\operatorname{cl}}
\newcommand{\hit}{\operatorname{hit}}

\usepackage{amsfonts}
\usepackage{amsmath}

\theoremstyle{theorem}
\newtheorem{lem}{Lemma}
\newtheorem{prop}[lem]{Proposition}
\newtheorem{cor}[lem]{Corollary}

\theoremstyle{definition}
\newtheorem{exam}[lem]{Example}

\newtheorem{defn}[lem]{Definition}

\begin{document}

\maketitle

\begin{abstract}
The bornological convergence structures that have been studied recently as generalizations of Attouch-Wets convergence define pretopologies on hyperspaces. In this paper we characterize the topological reflections of these pretopologies and translate the constructions necessary to define bornological convergence to a broader spectrum of hyperspace topologies.
\end{abstract}

\section{Introduction}

Attouch-Wets convergence is a convergence structure on the hyperspace of a metric or normed linear space with various applications in convex analysis. This type of convergence, also called bounded Hausdorff convergence, defines a topology that is coarser than the topology defined by the Hausdorff distance. The metrically bounded subsets play a fundamental role in this construction. Looking at the definition of Attouch-Wets convergence, the question arises whether this idea can be extended by using general bornologies instead of metrically bounded sets. In \cite{Lechicki} this generalisation of Attouch-Wets convergence was introduced and studied into detail. The construction that was given in this paper yields a neighbourhood filter for each element of the hyperspace, with one setback: this family of neighbourhood filters does in general not generate a topology, i.e. it is not certain that each neighbourhood filter has a base of open sets. Such an assignment of neighbourhood filters to the points of a set is called a \emph{pretopology}.

Determining wich conditions on a bornology are necessary and suffcient such that the corresponding pretopology is in fact topological is the main subject of \cite{Beer}. In this paper it is extensively explained what properties a bornology should have such that its bornological convergence structure (or its lower or upper part) is topological.

Within the category of pretopological spaces and continuous maps the topological spaces form a concretely reflective subcategory. This means that for each pretopology on a set $X$ there is a finest topology on $X$ such that a neighbourhood of a point in this topology is also a neighbourhood of that point in the pretopology. This brings us to the first question we want to answer in this paper: what is the topological reflection of a (lower or upper) bornological convergence structure as it was defined in \cite{Lechicki}?

The modifications that are necessary to transform a Hausdorff distance topology to Attouch-Wets convergence or (more general) bornological convergence can be performed on a much broader spectrum of hyperspace topologies. The second goal of this paper is to describe when such a \emph{bornological modification} of a hyperspace topology results in a new hyperspace topology (rather than a pretopology) and to give a description of the topological reflection in case the resulting structure is a pretopology.

\section{Preliminaries}

\subsection{Pretopological Spaces}

Pretopologies can be characterized in various ways. One way to characterize these objects is by a map from the set $X$ of points to the set of filters on $X$. The filter thus associated to a point will then be called the \emph{neighbourhood filter} of that point. A map between pretopological spaces is called \emph{continuous} iff the inverse image of a neighbourhood of $f(x)$ is a neighbourhood of $x$. A set $G\subseteq X$ is be called \emph{open} iff it is a neighbourhood of each of its points. Having assigned to each point a neighbourhood filter it is possible to define the closure operator $\cl$ as follows: $$a\in \cl(A) \ \Leftrightarrow \ \textrm{$X\setminus A$ is no neighbourhood of $X$}.$$ A closure operator that is constructed in this way has three essential properties.

\begin{itemize}
\item $A\subseteq \cl(A)$
\item $A\subseteq B \ \Rightarrow \ \cl(A)\subseteq \cl(B)$
\item $\cl(A\cup B)=\cl(A)\cup \cl(B)$
\end{itemize}

Conversely, each closure operator with these properties defines a pretopology by assigning to each point $x$ the following neighbourhood filter: $$\{V\subseteq X|x\not\in \cl(X\setminus V)\}.$$ This correspondence between pretopologies and closure operators is one-to-one. Using closure operators, continuous maps between pretopological spaces can be characterized as maps that send elements of $\cl(A)$ to elements of $\cl(f(A))$ for each $A\subseteq X$.

It is clear that each topology also defines a pretopology. The only difference between both being that in a topological space each neighbourhood filter has a base of open sets. In terms of closure operators we can say that a pretopology is a topology iff its closure operator is idempotent, i.e.

\begin{itemize}
\item $\cl(\cl(A))=\cl(A)$
\end{itemize}

Throughout this text the closure operator associated with a topology $T$ will be denoted as $\cl^T$. Although not every pretopology is a topology, we do have that the category of topological spaces and continuous maps is a concretely reflective subcategory of the category of pretopological spaces. For more information on this categorical terminology we refer the reader to \cite{Adamek} and \cite{Colebunders}. In a sense, we have that for each pretopology there is a topology that is the \textit{most similar} to it.  Concretely, it means that for each pretopological space $X$ there is a topological space $X^{\tau}$ with the same underlying set such that whenever $$f:X\to Y$$ is a continuous map and $Y$ a topological space, the map $$f:X^{\tau}\to Y$$ is also continuous. This topological reflection $X^{\tau}$ can be easily described using open sets. Its open sets are exactly the sets that are open in $X$. Characterizing the closure operator on $X^{\tau}$ is less straightforward and uses some cardinality arguments. If $\cl$ and $\cl'$ are closure operators, then we say that $\cl$ is finer than $\cl'$ iff $\cl(A)\subseteq \cl'(A)$ for each $A\subseteq X$. A closure operator is always finer than the closure operator of its topological reflection.

\subsection{Hyperspace Topologies}

For a collection $\mathcal{A}$ of subsets of a set $X$ we will define $\downarrow \mathcal{A}$ (resp. $\uparrow \mathcal{A}$) as $\{B|\exists A\in\mathcal{A}: B\subseteq A\}$ (resp. $\{B|\exists A\in \mathcal{A}:B\supseteq A\}$). If $\mathcal{A}$ is equal to $\downarrow \mathcal{A}$ (resp. $\uparrow \mathcal{A}$), then we will call $\mathcal{A}$ \emph{downward directed} (resp. \emph{upward directed}). The set $\downarrow \{A\}$ (resp. $\uparrow \{A\}$) will be simply denoted as $\downarrow A$ (resp. $\uparrow A$). A hyperspace topology in which all open sets are downward directed will be called an upper hyperspace topology. When all open sets are upward directed we will speak of a lower hyperspace topology. Many hyperspace topologies exist, each one with its specific applications and properties. Some well-known examples include the Wijsman topology, the Fell  topology and Vietoris topology (see for example \cite{BeerBook, Dimaio, Naimpally}).  Each one of these examples is the supremum of a lower and an upper hyperspace topology. 

The hyperspace topology that plays a key role in the seminal paper \cite{Lechicki} on bornological convergence is the topology defined by the Hausdorff distance. Each metric $d$ on a set $X$ defines an extended pseudometric $H_d$, called the \emph{Hausdorff distance}, on the hyperspace of $X$. The Hausdorff distance of two subsets $A$ en $B$ of the metric space $X$ is defined as the infemum of all $\epsilon>0$ for which $A\subseteq B^{\epsilon}$ and $B\subseteq A^{\epsilon}$. Here $A^{\epsilon}$ is the set of all points that have a $d$-distance less than $\epsilon$ from $A$. The topology defined by the Hausdorff distance is the supremum of a lower and an upper hyperspace topology. A base for the neighbourhoods of a set $A$ in the lower part are the sets of the form $\{B|A\subseteq B^{\epsilon}\}$. This topology will be denoted as $H^-$. The sets $\{B|B\subseteq A^{\epsilon}\}$ form a base of neighbourhoods for a set $A$ in the upper part, denoted as $H^+$.

A slight modification of the topology induced by the Hausdorff distance is the Attouch-Wets convergence \cite{Attouch}. In this topology a base for the neighbourhood filter of a set $A$ is given by the sets of the form $$\{B\subseteq X| A\cap S\subseteq B^{\epsilon} \textrm{ en } B\cap S\subseteq A^{\epsilon}\},$$ where $S$ is a subset of $X$ that is bounded for the metric $d$.

In \cite{Lechicki} the idea behind Attouch-Wets convergence was generalized some more, in the sense that the authors no longer restricted themselves to the class of metrically bounded sets, but considered an arbitrary order-theoretic ideal $\mathcal{S}$ of subsets of $X$. Throughout this paper we will often refer to the elements of $\mathcal{S}$ as the \emph{bounded} sets. An essential difference between Attouch-Wets convergence and this newly defined \emph{bornological convergence} is that the latter need not be a topology, but is in general a pretopology. In \cite{Lechicki} and in later papers (see \cite{Beer}) on bornological convergence the study of this type of convergence is split up in three parts. Lower bornological convergence is defined by the closure operator $\cl_{\mathcal{S}}^{H^-}$. A set $\cl_{\mathcal{S}}^{H^-}(\mathcal{A})$ contains all $A\subseteq X$ that have the property that for each $\epsilon>0$ and $\mathcal{S}$ there is a $B\in\mathcal{A}$ for which $$A\cap S\subseteq B^{\epsilon}.$$ The closure operator $\cl_{\mathcal{S}}^{H^+}$ defines the notion of upper bornological convergence. A set $A$ is contained in the upper bornological closure of a set $\mathcal{A}$ iff for each $\epsilon>0$ and $S\in\mathcal{S}$ there is a $B\in\mathcal{A}$ such that $$B\cap S\subseteq A^{\epsilon}.$$ Bornological convergence finally, is defined by the closure operator $\cl_{\mathcal{S}}^{H}$. The bornological closure of a set $\mathcal{A}$ contains all sets $A$ that satisfy the property that for each $\epsilon>0$ and $S\in \mathcal{S}$ there is a $B\in\mathcal{A}$ such that $$A\cap S\subseteq B^{\epsilon} \textrm{ and } B\cap S\subseteq A^{\epsilon}.$$ Properties of these types of convergence were studied in \cite{Rodriguez}.

\section{Lower Bornological Convergence}

Throughout this text $(X,d)$ will be a metric space and $\mathcal{S}$ an ideal of subsets of $X$. As mentioned before lower bornological convergence is defined by the closure operator $\cl_{\mathcal{S}}^{H^-}$. The topological reflection of this pretopology will be denoted as $H^-(\mathcal{S})$. We will begin by giving a description of the open sets in $H^-(\mathcal{S})$. After that we will take a closer look at how this construction can be translated to a more general setting of lower hyperspace topologies.

\begin{defn}
We say that a collection of sets $\mathcal{S}'$ satisfies condition $(\clubsuit)$ iff for each $S\in\mathcal{S}'$ there is an $\epsilon>0$ such that whenever $S\subseteq B^{\epsilon}$ there is an $S_B\in\mathcal{S}'$ that is a subset of $B$.
\end{defn}

\begin{prop}
The open sets in $H^-(\mathcal{S})$ are exactly the sets $$\bigcup_{S\in\mathcal{S}'}\uparrow S$$ where $\mathcal{S}'$ is a subset of $\mathcal{S}$ that satisfies condition $(\clubsuit)$.
\begin{proof}
Let $\mathcal{S}'$ be a subset of $\mathcal{S}$ that satisfies condition $(\clubsuit)$ and $A$ a set that contains an $S_0\in\mathcal{S}'$. By assumption we can find an $\epsilon>0$ such that whenever $S_0\subseteq B^{\epsilon}$ we can find an $S_B\subseteq B$ such that $S_B\in\mathcal{S}'$. This yields that $\{B|A\cap S_0\subseteq B^{\epsilon}\}$ is a subset of $$\bigcup_{S\in\mathcal{S}'}\uparrow S$$ and that the latter is open in $H^-(\mathcal{S})$.

Now let $\mathcal{G}$ denote an open set in $H^-(\mathcal{S})$ and define $\mathcal{S}'$ as the set of all bounded elements of $\mathcal{G}$. For each $A\in\mathcal{G}$ we can find an $S\in\mathcal{S}$ and an $\epsilon>0$ such that $\{B|A\cap S\subseteq B^{\epsilon}\}$. This means that $A\cap S$ is a bounded set that is a subset of $A$ and an element of $\mathcal{G}$. Hence we obtain that $\mathcal{G}$ is equal to $$\bigcup_{S\in\mathcal{S}'}\uparrow S.$$ Take an $S_0$ in $\mathcal{S}'$ and choose an $\epsilon>0$ such that $\{B|S\subseteq B^{\epsilon}\}\subseteq \mathcal{G}$. For each $B$ that satisfies $S\subseteq B^{\epsilon}$ we have that $B$ is an element of $\mathcal{G}$ and therefore contains an $S_B\in\mathcal{S}'$. This implies that $\mathcal{S}'$ satisfies condition $(\clubsuit)$.
\end{proof}
\end{prop}

The topology $H^-(\mathcal{S})$ itself defines a closure operator. The next question we want to answer is whether we can find an ideal $\mathcal{S}^-$ of subsets of $X$ such that the closure operator defined by $H^-(\mathcal{S})$ is equal to $\cl_{\mathcal{S}^-}^{H^-}$. If such an ideal exists, then by assumption the closure operator $\cl_{\mathcal{S}^-}^{H^-}$ is topological. It was established in \cite{Beer} that this is equivalent to saying that for each $S\in\mathcal{S}^-$ we have $$\forall \epsilon>0\exists\delta>0\forall S\subseteq A^{\delta}\exists S_A\in\mathcal{S}^-:\textrm{$S_A\subseteq A$ and $S\subseteq S_A^{\epsilon}$}.$$

\begin{prop}
Each $\mathcal{S}$ contains a largest $\mathcal{S}^*$ for which $H^-(\mathcal{S}^*)$ is topological.
\begin{proof}
For an arbitrary ideal $\mathcal{S}'$ we will define $\widehat{\mathcal{S}'}$ as the set of all $S\in\mathcal{S}'$ that satisfy $$\forall \epsilon>0\exists\delta>0\forall  A^{\delta}\supseteq S\exists S_A\in\mathcal{S}':\textrm{$S_A\subseteq A$ and $S\subseteq S_A^{\epsilon}$}.$$ Define $\mathcal{S}_0$ as $\mathcal{S}$ and $\mathcal{S}_{\gamma}$ as $\widehat{\mathcal{S}_{\alpha}}$ if $\gamma$ is equal to a successor ordinal $\alpha+1$ and as $$\bigcap_{\beta<\gamma}\mathcal{S}_{\beta}$$ if $\gamma$ is a limit ordinal. There exists an ordinal $\gamma$ for which $\mathcal{S}_{\gamma}=\mathcal{S}_{\gamma+1}$. Denote the latter as $\mathcal{S}^*$. By definition we have that $\mathcal{S}^*$ satisfies the necessary conditions for which $H^-(\mathcal{S}^*)$ is topological.

Let $\mathcal{S}'$ be an ideal that is contained in $\mathcal{S}$ and for which $H^-(\mathcal{S}')$ is topological. This yields that $\mathcal{S}'$ is a subset of each $\mathcal{S}_{\beta}$ for all $\beta<\gamma$. Hence $\mathcal{S}'$ is a subset of $\mathcal{S}_{\gamma}$ and thus we obtain that $\mathcal{S}'$ is a subset of $\mathcal{S}^*$.
\end{proof}
\end{prop}

\begin{defn}
$\mathcal{S}^{tb}$ is defined as the set of all $S'\subseteq X$ that satisfy the following condition: $$\forall \epsilon>0 \ \exists S\in \mathcal{S}: S'\subseteq S^{\epsilon}.$$ 
\end{defn}

\begin{prop}
A closure operator $\cl_{\mathcal{S}}^{H^-}$ is finer than a closure operator $\cl_{\mathcal{S}'}^{H^-}$  iff $\mathcal{S}'$ is contained in $\mathcal{S}^{tb}$.
\begin{proof}
See \cite{Lechicki}.
\end{proof}
\end{prop}

Because $(\mathcal{S}^{tb})^{tb}$ is equal to $\mathcal{S}^{tb}$ this also implies that bornological convergences defined by $\mathcal{S}$ and $\mathcal{S}^{tb}$ are equal.

\begin{defn}
For an arbitrary ideal $\mathcal{S}$ we will denote $(\mathcal{S}^{tb})^*$ as $\mathcal{S}^-$.
\end{defn}

\begin{prop}
If the closure operator associated with the topology $H^-(\mathcal{S})$ is equal to $\cl^{H^-}_{\mathcal{S}'}$ for some bornology $\mathcal{S}'$, then it is equal to $\cl^{H^-}_{\mathcal{S}^-}$.
\begin{proof}
Because $\cl^{H^-}_{\mathcal{S}^-}$ is by definition topological and $H^-(\mathcal{S})$ is the topological reflection of the pretopological closure operator $\cl^{H^-}_{\mathcal{S}}$ we obtain that $\cl_{\mathcal{S}'}^{H^-}$ is finer than $\cl_{\mathcal{S}^-}^{H^-}$ and that $\mathcal{S}^-\subseteq (\mathcal{S}')^{tb}$. On the other hand, we have that $\cl_{\mathcal{S}}^{H^-}$ is finer than $\cl_{\mathcal{S}'}^{H^-}$ and that $(\mathcal{S}')^{tb}\subseteq \mathcal{S}^{tb}$. Because of the previous proposition this yields $(\mathcal{S}')^{tb}\subseteq \mathcal{S}^-$. We can now conclude that the closure operators $cl^{H^-}_{\mathcal{S}'}$ and $cl^{H^-}_{\mathcal{S}^-}$ are equal and that the latter is equal to the closure operator associated with the topology $H^-(\mathcal{S})$.
\end{proof}
\end{prop}

The following example shows that $\cl_{\mathcal{S}^-}^{H^-}$ does not need to be equal to the closure operator defined by $H^-(\mathcal{S})$.

\begin{exam}
We will work in $\mathbb{C}$ with the Euclidian metric and define $B_{\epsilon}$ as $$\{z\in\mathbb{C}|-\epsilon<\Im(z)<\epsilon\}.$$ Let $\mathcal{S}$ be the ideal defined as $$\{S\subseteq \mathbb{C}|\forall \epsilon>0: \textrm{$S\setminus B_{\epsilon}$ is bounded}\}.$$ It is easy to verify that $\mathcal{S}^{tb}$ is equal to $\mathcal{S}$ itself. Let $S$ be an unbounded element of $\mathcal{S}$. For each $\delta>0$ it holds that $S\subseteq A^{\delta}$, where $A$ is defined as $$S+\frac{i\delta}{2}.$$ There is, however, no element $S_A$ in $\mathcal{S}$ that is contained in $A$ and that satisfies $S\subseteq (S_A)^1$. This means that all elements in $\mathcal{S}^-$ are bounded. Since $\cl_{\mathcal{S}_b}^{H^-}$ is equal to Attouch-Wets convergence on the subsets of $\mathbb{C}$, where $\mathcal{S}_b$ is the set of all bounded subsets of $\mathbb{C}$, and this type of convergence is always topological, we obtain that $\mathcal{S}^-$ is equal to $\mathcal{S}_b$.

Define $\mathcal{G}$ as the set of all subsets of $\mathbb{C}$ except the finite subsets of $\mathbb{Z}$. This set is open in $H^-(\mathcal{S})$. For each element $A\in\mathcal{G}$ that has an infinite intersection with $\mathbb{Z}$ we have that $\{B|A\cap \mathbb{Z}\subseteq B^{\frac{1}{2}}\}$ is a neighbourhood of $A$ that is contained in $\mathcal{G}$. If $A$ has a finite intersection with $\mathbb{Z}$, then the same holds for $\{B|A\cap S\subseteq B^{\delta}\}$, where $S$ is an arbitrary finite, non-empty subset of $A\setminus \mathbb{Z}$ and $\delta$ is chosen smaller than the distance between $S$ and $\mathbb{Z}$. The set $\mathcal{G}$ is, however, not open in $H^-(\mathcal{S}_b)$. Indeed, there is no $S\in\mathcal{S}_b$ and no $\epsilon>0$ such that $\{B|\mathbb{Z}\cap S\subseteq B^{\epsilon}\}\subseteq \mathcal{G}$, while $\mathbb{Z}$ is an element of $\mathcal{G}$.
\end{exam}

Looking at the definition of lower bornological convergence we see that a set $A$ is in the lower bornological closure of $\mathcal{A}$ iff each element of $\mathcal{S}$ that is contained in $A$ is an element of the lower Hausdorff closure of $\mathcal{A}$. In what follows we will generalize this idea to arbitrary lower hyperspace topologies and study the properties of the pretopologies that can be defined this way.

\begin{defn}
For each $A\subseteq X$ we will denote the set of all elements in $\mathcal{S}$ that are contained in $A$ as $\mathcal{S}(A)$. For a lower hyperspace topology $T$ and a set $\mathcal{A}$ we define $\cl_{\mathcal{S}}^T(\mathcal{A})$ as the set that contains all $A\subseteq X$ that satisfy $$S\in\cl^T(\mathcal{A})$$ whenever $S\in\mathcal{S}(A)$.
\end{defn}

\begin{prop}
$\cl^T_{\mathcal{S}}$ is a closure operator.
\begin{proof}
Let $\mathcal{A}$ be a collection of subsets of $X$. Take $A\in\mathcal{A}$ and $S\in\mathcal{S}(A)$. Because $T$ is a lower hyperspace topology we have that each neighbourhood of $S$ contains $A$ and therefore $A\in \cl^T_{\mathcal{S}}(\mathcal{A})$. It can be easily seen that $\cl^T_{\mathcal{S}}$ preserves inclusions. Suppose that $A$ is not contained in $\cl^T_{\mathcal{S}}(\mathcal{A})\cup \cl^T_{\mathcal{S}}(\mathcal{B})$. By definition we can find $S_1,S_2\in\mathcal{S}(A)$ such that $S_1\not\in \cl^T(\mathcal{A})$ and $S_2\not\in \cl^T(\mathcal{B})$. This yields that $S_1\cup S_2$ is an element of $\mathcal{S}(A)$ that is not contained in $\cl^T(\mathcal{A}\cup\mathcal{B})$. Hence we obtain that $A$ is no element of $\cl^T(\mathcal{A}\cup\mathcal{B})$ and that $\cl^T_{\mathcal{S}}(\mathcal{A})\cup\cl^T_{\mathcal{S}}(\mathcal{B})$ and $\cl^T_{\mathcal{S}}(\mathcal{A}\cup \mathcal{B})$ are equal.
\end{proof}
\end{prop}

\begin{prop}
$\cl^T_{\mathcal{S}}$ is finer than $\cl^T_{\mathcal{S}^*}$ iff each $S^*\in\mathcal{S}^*$ is an element of $\cl^T(\mathcal{S}(S^*))$.
\begin{proof}
Suppose that $\cl^T_{\mathcal{S}}$ is finer than $\cl^T_{\mathcal{S}^*}$ and take $S^*\in \mathcal{S}^*$. We automatically have $S^*\in\cl^T_{\mathcal{S}}(\mathcal{S}(S^*))$ and, by assumption, this yields $S^*\in cl^T_{\mathcal{S}^*}(\mathcal{S}(S^*))$. Since $S^*\in\mathcal{S}^*$ the latter is equivalent to $S^*\in\cl^T(\mathcal{S}(S^*))$.

The other way around, assume that each $S^*\in\mathcal{S}^*$ is contained in $\cl^T(\mathcal{S}(S^*))$ and that $A$ is an element of $\cl^T_{\mathcal{S}}(\mathcal{A})$. Because $T$ is a lower hyperspace topology we get $S^*\in\cl^T_{\mathcal{S}}(\mathcal{A})$ for each $S^*\in\mathcal{S}^*(A)$. The latter yields that $\mathcal{S}(S^*)$ is a subset of $\cl^T(\mathcal{A})$. Since $S^*\in\cl^T(\mathcal{S}(S^*))$ we obtain $S^*\in \cl^T(\mathcal{A})$. This means that $A$ is an element of $\cl^T_{\mathcal{S}^*}(\mathcal{A})$.
\end{proof}
\end{prop}

\begin{defn}
For any two sets $A,B\subseteq X$ we define $[A,B]$ as the collection of all sets $C\subseteq X$ that contain $A$ and are disjoint from $B$. The sets of the form $[S,S']$, with $S,S'\in\mathcal{S}$, are a base for a topology. This topology will be denoted $\tau(\mathcal{S})$ and $\cl_{\mathcal{S}}$ denotes the associated closure operator.
\end{defn}

The topology we just defined will play an important role in both the description of the topological reflection of the lower and the upper bornological convergence structures.

\begin{prop}
A set $A$ is an element of $\cl_{\mathcal{S}}(\mathcal{B})$ iff for each $S\in\mathcal{S}$ there is a $B\in\mathcal{B}$ such that $A\cap S=B\cap S$.
\begin{proof}
Suppose $A\in\cl_{\mathcal{S}}(\mathcal{B})$ and take $S\in\mathcal{S}$. Define $S_1$ as $A\cap S$ and $S_2$ as $S\setminus A$. By definition we can find a $B\in\mathcal{B}$ such that $B\in [S_1,S_2]$, This yields $A\cap S=B\cap S$.

To prove the converse, assume that for each $S\in\mathcal{S}$ there is a $B\in\mathcal{B}$ such that $A\cap S=B\cap S$. Let $A$ be an element of $[S,S']$ and define $S^*$ as $S\cup S'$. We can now find a $B\in\mathcal{B}$ such that $A\cap S^*=B\cap S^*$
\end{proof}
\end{prop}

We give the following lemma without proof.

\begin{lem}
Let $\mathcal{B}$ be a collection of subsets of $X$.
\begin{enumerate}
\item If $\mathcal{B}$ is downward directed, then $A\in\cl_{\mathcal{S}}(\mathcal{B})$ iff for each $S\in\mathcal{S}$ there is a $B\in\mathcal{B}$ such that $A\cap S\subseteq B$.
\item If $\mathcal{B}$ is upward directed, then $A\in\cl_{\mathcal{S}}(\mathcal{B})$ iff for each $S\in\mathcal{S}$ there is a $B\in\mathcal{B}$ such that $B\subseteq A\cap S$.
\end{enumerate}
\end{lem}

\begin{prop}
The closure operator $\cl_{\mathcal{S}}^T$ is equal to $\cl_{\mathcal{S}}\circ \cl^T$.
\begin{proof}
Let $A$ be an element of $\cl_{\mathcal{S}}^T(\mathcal{A})$ and take $S\in\mathcal{S}$. By definition we have $A\cap S\in \cl^T(\mathcal{A})$. Hence we obtain that $A$ is an element of $\cl_{\mathcal{S}}(\cl^T(\mathcal{A}))$.

If $A$ is an element of $\cl_{\mathcal{S}}(\cl^T(\mathcal{A}))$ and we take $S\in\mathcal{S}(A)$, then there is a $B\in \cl^T(\mathcal{A})$ such that $S\subseteq B$. This yields that $S$ is itself an element of $\cl^T(\mathcal{A})$ and thus we get that $A$ is contained in $\cl_{\mathcal{S}}^T(\mathcal{A})$.
\end{proof}
\end{prop}

\begin{cor}
$\cl_{\mathcal{S}}^{H^-}$ is equal to $\cl_{\mathcal{S}}\circ \cl^{H^-}$.
\end{cor}

\begin{prop}
The closure operator $\cl_{\mathcal{S}}^T$ is topological iff $\cl_{\mathcal{S}}(\mathcal{A})$ is closed in $T$ whenever $\mathcal{A}$ is closed in $T$.
\begin{proof}
If $\cl_{\mathcal{S}}^T$ is topological and $\mathcal{A}$ is closed in $T$, then the following inclusions hold:
\begin{eqnarray}
\cl^T(\cl_{\mathcal{S}}(\mathcal{A})) & = & \cl^T(\cl_{\mathcal{S}}(\cl^T(\mathcal{A}))) \nonumber \\
& \subseteq & \cl_{\mathcal{S}}^T(\cl_{\mathcal{S}}^T(\mathcal{A})) \nonumber \\
& = & \cl_{\mathcal{S}}^T(\mathcal{A}) \nonumber \\
& = & \cl_{\mathcal{S}}(\cl^T(\mathcal{A})) \nonumber \\
& = & \cl_{\mathcal{S}}(\mathcal{A}) \nonumber
\end{eqnarray}
This yields that $\cl_{\mathcal{S}}(\mathcal{A})$ is closed in $T$.

To prove the necessity of this condition we notice the following:
\begin{eqnarray}
\cl_{\mathcal{S}}^T(\cl_{\mathcal{S}}^T(\mathcal{A})) & = & \cl_{\mathcal{S}}(\cl^T(\cl_{\mathcal{S}}(\cl^T(\mathcal{A})))) \nonumber \\
& = & \cl_{\mathcal{S}}(\cl_{\mathcal{S}}(\cl^T(\mathcal{A}))) \nonumber \\
& = & \cl_{\mathcal{S}}(\cl^T(\mathcal{A})) \nonumber 
\end{eqnarray}
\end{proof}
\end{prop}

\begin{prop}
The topology $T\wedge \tau(\mathcal{S})$ is the topological reflection of the pretopology defined by $\cl_{\mathcal{S}}^T$.
\begin{proof}
It is, by definition, clear that all closed sets in $T\wedge \tau(\mathcal{S})$ are closed for the closure operator $\cl_{\mathcal{S}}^T$. Conversely, let $\mathcal{A}$ be a set that is closed for this closure operator. First of all, we have the following inclusions: 
\begin{eqnarray}
\cl_{\mathcal{S}}(\mathcal{A}) & \subseteq & \cl_{\mathcal{S}}(\cl^T(\mathcal{A})) \nonumber \\
& = & \cl_{\mathcal{S}}^T(\mathcal{A}) \nonumber \\
& = & \mathcal{A} \nonumber
\end{eqnarray}
Furthermore, we have:
\begin{eqnarray}
\cl^T(\mathcal{A}) & \subseteq & \cl_{\mathcal{S}}(\cl^T(\mathcal{A})) \nonumber \\
& = & \cl_{\mathcal{S}}^T(\mathcal{A}) \nonumber \\
& = & \mathcal{A} \nonumber
\end{eqnarray}
Hence we obtain that $\mathcal{A}$ is closed in $T\wedge \tau(\mathcal{S})$.
\end{proof}
\end{prop}

\begin{cor}
$H^-(\mathcal{S})$ is equal to $H^-\wedge \tau(\mathcal{S})$.
\end{cor}

\begin{defn}
From here on we will denote the topology $T\wedge \tau(\mathcal{S})$ as $T(\mathcal{S})$.
\end{defn}

\begin{defn}
The topology $T(\mathcal{S})$ satisfies the property that a set $A$ is in the closure of $\mathcal{A}$ iff all of its bounded subsets belong to the closure of $\mathcal{A}$. We will call such hyperspace topologies \emph{boundedly generated}.
\end{defn}

\begin{prop}
$T(\mathcal{S})$ is the finest boundedly generated topology coarser than $T$.
\begin{proof}
Let $\mathcal{A}$ be a closed set in such a hyperspace topology and take $A\in \cl_{\mathcal{S}}(\mathcal{A})$. Since we can find a $B\in\mathcal{A}$ such that $S\subseteq B$ for each $S\in\mathcal{S}(A)$ and $\mathcal{A}$ is closed in a lower hyperspace topology we obtain that $S$ is an element of $\mathcal{A}$ for each $S\in\mathcal{S}$. This of course yields that each bounded subset of $A$ is contained in $\mathcal{A}$ and thus that $A$ is an element of $\mathcal{A}$. Hence we obtain that $\mathcal{A}$ is closed in $T(\mathcal{S})$.
\end{proof}
\end{prop}

\section{Upper Bornological Convergence}

Similar to what we did in the previous section with lower bornological convergence, we will now give a characterization of the open sets in the topological reflection $H^+(\mathcal{S})$ of the pretopology defined by $\cl_{\mathcal{S}}^{H^+}$. 

\begin{defn}
An ideal $\mathcal{S}$ will be called \emph{stable under small enlargements} iff for each $S\in\mathcal{S}$ we can find an $\epsilon>0$ such that $S^{\epsilon}\in\mathcal{S}$.

A set $C$ will be called \emph{cobounded} iff its complement is bounded.
\end{defn}

\begin{prop}
The open sets in $H^+(\mathcal{S})$ are exactly the sets $$\bigcup_{C\in\mathcal{C}}\downarrow C$$ where $\mathcal{C}$ is a collection of cobounded sets that is stable under small enlargements.
\begin{proof}
Let $C$ be a collection of cobounded sets that is stable under small enlargements. Suppose that $A$ is an element in $\downarrow C$ for some $C\in\mathcal{C}$ and that $C^{\epsilon}\in\mathcal{C}$. If we define $S$ as $X\setminus C^{\epsilon}$, then $\{B|B\cap S\subseteq A^{\epsilon}\}$ is a subset of $\downarrow(C^{\epsilon})$ and therefore contained in $$\bigcup_{C\in\mathcal{C}}\downarrow C.$$

Take an open set $\mathcal{G}$ in $H^+(\mathcal{S})$ and let $A$ be an element of $\mathcal{G}$. Because $\mathcal{G}$ is open we can find an $S\in\mathcal{S}$ and an $\epsilon>0$ such that $\{B|B\cap S\subseteq A^{\epsilon}\}$ is contained in $\mathcal{G}$. Define $C$ as $X\setminus \left(S\setminus A^{\epsilon}\right)$. The set $\downarrow C$ contains $A$ and is a subset of $\mathcal{G}$. Now let $\mathcal{C}$ be the collection of all cobounded sets in $\mathcal{G}$. For each $C\in\mathcal{C}$ we can find a $\delta>0$ such that $\{B|B\subseteq C^{\delta}\}\subseteq \mathcal{G}$. Hence we obtain that $C^{\delta}$ is itself a cobounded element of $\mathcal{G}$. This yields that $\mathcal{C}$ is stable under small enlargements. Moreover, we can conclude that $\mathcal{G}$ is equal to the union of all sets $\downarrow C$, with $C\in\mathcal{C}$.
\end{proof}
\end{prop}

\begin{defn}
A hyperspace topology that has a base that consists of sets of the form $\downarrow G$ will be called a \emph{miss topology}. This terminology is chosen because of the fact that all elements in $\downarrow G$ miss the complement of $G$. Examples can be found in \cite{Naimpally}. It is clear that each miss topology is an upper hyperspace topology.
\end{defn}

\begin{prop}
$H^+(\mathcal{S})$ is the finest miss topology that is contained in $H^+$ and for which the cobounded sets are dense.
\begin{proof}
Let $T$ be a miss topology that is contained in $H^+$ and for which the cobounded sets are dense. Take an open set $\mathcal{G}$ in this topology and an $A$ in $\mathcal{G}$. By assumption we can now find a $G_A$ such that $\downarrow G_A$ is a neighbourhood of $A$ that is contained in $\mathcal{G}$. Because the cobounded sets are dense we obtain that $\downarrow G_A$ contains a cobounded set. This implies that $G_A$ is itself cobounded. If we now denote the collection of all cobounded sets in $\mathcal{G}$ as $\mathcal{C}$ we obtain that $\mathcal{G}$ is equal to $$\bigcup_{C\in\mathcal{C}}\downarrow C.$$ For each $C\in\mathcal{C}$ we can find an $\epsilon>0$ such that $\{B|B\subseteq C^{\epsilon}\}\subseteq \mathcal{G}$ and thus we have that $C^{\epsilon}$ is an enlargement of $C$ that is cobounded and contained in $G$. This yields that $\mathcal{C}$ is closed under small enlargements and that $G$ is an open set in the topological reflection of $H^+(\mathcal{S})$.
\end{proof}
\end{prop}

Contrary to the situation with lower bornological convergence, the closure operator defined by the topological reflection of an upper bornological convergence structure is always equal to an upper bornological closure operator. It was established in \cite{Beer} that an upper bornological closure operator $\cl_{\mathcal{S}}^{H^+}$ is topological iff $\mathcal{S}$ is stable under small enlargements.

\begin{defn}
For an ideal $\mathcal{S}$ we define $\mathcal{S}^+$ as the set $\{S\in\mathcal{S}|\exists \epsilon>0:S^{\epsilon}\in\mathcal{S}\}$. This set is again an ideal.
\end{defn}

\begin{prop}
$\mathcal{S}^+$ is stable under small enlargements.
\begin{proof}
Take $S\in\mathcal{S}^+$ and $\epsilon>0$ such that $S^{\epsilon}\in\mathcal{S}$. Since $(S^{\frac{\epsilon}{2}})^{\frac{\epsilon}{2}}$ is contained in $S^{\epsilon}$ we get that the former is an element of $\mathcal{S}$. This means that $S^{\frac{\epsilon}{2}}$ is an enlargement of $S$ that is an element of $\mathcal{S}^+$.
\end{proof}
\end{prop}

\begin{prop}
The closure operator associated with $H^+(\mathcal{S})$ is equal to $\cl_{\mathcal{S}^+}^{H^+}$.
\begin{proof}
The previous proposition yields that the closure operator associated with $H^+(\mathcal{S})$ is finer than $\cl_{\mathcal{S}}^{H^+}$. An open set $\mathcal{G}$ in $H^+(\mathcal{S})$ is of the form $$\bigcup_{C\in\mathcal{C}}\downarrow C$$ where $\mathcal{C}$ is a collection of cobounded sets that is stable under small enlargements. Take $A\in\mathcal{G}$ and $C\in\mathcal{C}$ such that $A\subseteq C$. By assumption there is an $\epsilon>0$ such that $C^{\epsilon}\in\mathcal{C}$. If we define $S$ as $X\setminus C^{\epsilon}$, then $S$ is a bounded set for which $S^{\epsilon}$ is again bounded. This means that $S\in\mathcal{S}^+$. The set $\{B|B\cap S\subseteq A^{\epsilon}\}$, being equal to $\downarrow(C^{\epsilon})$, is a $H^+(\mathcal{S}^+)$-neighbourhood of $A$ that is contained in $\mathcal{G}$. Because $A$ was arbitrary we obtain that $\mathcal{G}$ is open in $H^+(\mathcal{S}^+)$.
\end{proof}
\end{prop}

The following results describe how the construction of an upper bornological convergence structure can be generalized to the setting of arbitrary upper hyperspace topologies. To avoid an abundance of new notation we will use $\cl_{\mathcal{S}}^T$ again to describe the bornological modification of an upper hyperspace topology. Since the indiscrete topology is the only upper hyperspace topology that is also a lower hyperspace topology and both definitions coincide for this particular topology there will be no room for ambiguity.

\begin{defn}
For an upper hyperspace topology $T$ we define $\cl_{\mathcal{S}}^T(\mathcal{A})$ as follows: $$A\in \cl_{\mathcal{S}}^T(\mathcal{A})\Leftrightarrow \forall S\in\mathcal{S}: A\in \cl^T(\mathcal{A}|_S)$$ where $\mathcal{A}|_S$ is defined as $\{B\cap S|B\in\mathcal{A}\}$.
\end{defn}

\begin{prop}
$\cl_{\mathcal{S}}^T$ is a closure operator.
\begin{proof}
Because $T$ is an upper hyperspace topology we obtain that each $\mathcal{A}$ is a subset of $\cl_{\mathcal{S}}^T(\mathcal{A})$. That $\cl_{\mathcal{S}}^T$ is order preserving is clear. Suppose that $A$ is no element of $\cl_{\mathcal{S}}^T(\mathcal{A})$ or $\cl_{\mathcal{S}}^T(\mathcal{A}')$. By definition this yields that we can find $S$ and $S'$ in $\mathcal{S}$ such that $A\not\in \cl^T(\mathcal{A}|_S)$ and $A\not\in \cl^T(\mathcal{A}'|_{S'})$. Since $\cl^T((\mathcal{A}\cup \mathcal{A}')|_{S\cup S'})$ is a subset of $\cl^T(\mathcal{A}|_S)\cup \cl^T(\mathcal{A}'|_{S'})$ we obtain that $A$ is no element of $\cl_{\mathcal{S}}^T(\mathcal{A}\cup \mathcal{A}')$.
\end{proof}
\end{prop}

\begin{prop}
If $T$ is a miss topology, then for each $\mathcal{B}$ the set $\cl_{\mathcal{S}}^T(\mathcal{B})$ is contained in $\cl_{\mathcal{S}}(\cl^T(\mathcal{B}))$.
\begin{proof}
Let $A$ be an element of the former set and take $S\in\mathcal{S}$. Define $A'$ as $A\cup (X\setminus S)$. For each $T$-neighbourhood $\downarrow G$ of $A'$ we can find a $B\in\mathcal{B}$ such that $B\cap S\subseteq G$. This implies that $A'$ is an element of $\cl^T(\mathcal{B})$. Since $A'\cap S\subseteq A$ we have that $A$ is in fact an element of $\cl_{\mathcal{S}}(\cl^T(\mathcal{B}))$.
\end{proof}
\end{prop}

\begin{prop}
For each $\mathcal{B}$ the set $\cl_{\mathcal{S}}^T(\mathcal{B})$ contains $\cl^T(\cl_{\mathcal{S}}(\mathcal{B}))$.
\begin{proof}
Let $A$ be an element of the latter set, $\downarrow G$ a $T$-neighbourhood of $A$ and $S\in\mathcal{S}$. By definition we can find a $B'\in\cl^+_{\mathcal{S}}(\mathcal{B})$ such that $B'\subseteq G$. Because we can find a $B\in\mathcal{B}$ such that $B\cap S\subseteq B'$, and thus $B\cap S\subseteq G$, we obtain that $A$ is an element of $\cl_{\mathcal{S}}^T(\mathcal{B})$.
\end{proof}
\end{prop}

\begin{prop}
$\cl^T_{\mathcal{S}}$ is finer than $\cl^T_{\mathcal{S}^*}$ iff the empty set is an element of $$\cl^T\left(\{S^*\setminus S|S\in\mathcal{S}\}\right)$$ for each $S^*\in\mathcal{S}^*$.
\begin{proof}
Let $\mathcal{C}$ be the collection of all $\mathcal{S}$-cobounded sets. By definition we have that the empty set is an element of $\cl^T(\mathcal{C}|_S)$ whenever $S\in\mathcal{S}$. This yields that if $\cl^T_{\mathcal{S}}$ is finer than $\cl^T_{\mathcal{S}^*}$, the empty set is an element of $\cl^T(\mathcal{C}|_{S^*})$ for each $S^*\in\mathcal{S}^*$. Because $\mathcal{C}|_{S^*}$ is equal to $\{S^*\setminus S|S\in\mathcal{S}\}$ we obtain that the given condition is necessary.

To prove its sufficiency assume that $A$ is an element of $\cl^T(\mathcal{A}|_S)$ for each $S\in\mathcal{S}$ and let $S^*$ be an element of $\mathcal{S}^*$. Because $T$ is an upper hyperspace topology we have that $A$ is an element of $\cl^T\left(\{S^*\setminus S|S\in\mathcal{S}\}\right)$. This means that for a neighbourhood $\downarrow V$ of $A$ we can find an $S\in\mathcal{S}$ such that $S^*\setminus S\subseteq V$. Take an $A'\in\mathcal{A}$ such that $A'\cap S\subseteq V$. We then obtain that $A'\cap S^*$ is contained in $V$. Hence we can conclude that $A$ is an element of $\cl^T(\mathcal{A}|_{S^*})$.
\end{proof}
\end{prop}

\begin{prop}
If $T$ is a miss hyperspace topology, then $T\wedge \tau(\mathcal{S})$ is the topological reflection of the closure operator $\cl_{\mathcal{S}}^T$.
\begin{proof}
Let $\mathcal{A}$ be closed for the closure operator $\cl_{\mathcal{S}}^T$. Suppose that each set $\mathcal{G}$ that is open for the topology $T$ and contains $A$ intersects with $\mathcal{A}$. Let $B$ be an element of this intersection. Because $T$ is an upper hyperspace topolog the set $B\cap S$ is an element of $\mathcal{G}$ for each $S\in\mathcal{S}$. This yields that $A$ is contained in $\cl_{\mathcal{S}}^T(\mathcal{A})$ and thus in $\mathcal{A}$. Hence we obtain that $\mathcal{A}$ is closed in $T$. Now assume $A\in\cl_{\mathcal{S}}(\mathcal{A})$. This means that for each $S\in\mathcal{S}$ we can find a $B\in\mathcal{A}$ such that $A\cap S=B\cap S$. Because $T$ is an upper hyperspace topology this yields $A\in\cl^T(\mathcal{A}|_S)$, so we can conclude that $\mathcal{A}$ is closed in $T\wedge \tau(\mathcal{S})$.

Conversely, let $\mathcal{A}$ be a closed in $T\wedge\tau(\mathcal{S})$ and take $A\in\cl_{\mathcal{S}}^T(\mathcal{A})$. If $\mathcal{G}$ is an open set in $T\wedge \tau(\mathcal{S})$ that contains $A$, then $\mathcal{G}$ contains a cobounded set $C$ that contains $A$. Since $T$ is a miss hyperspace topology we know that we can find a $T$-neighbourhood of $C$ of the form $\downarrow C^{\ast}$ that is contained in $\mathcal{G}$. By assumption the sets $\downarrow C^{\ast}$ and $\mathcal{A_{S^{\ast}}}$, where $S^{\ast}$ is the complement of $C^{\ast}$, must have a non-empty intersection. This implies that $\mathcal{A}$ contains a set that is disjoint from $S^{\ast}$ and therefore an element of $\mathcal{G}$. This leaves us to conclude that $A$ is an element of $\mathcal{A}$ and that the latter is closed for $\cl_{\mathcal{S}}^T$.
\end{proof}
\end{prop}

\begin{cor}
$H^+(\mathcal{S})$ is equal to $H^+\wedge \tau(\mathcal{S})$.
\end{cor}

\begin{defn}
Like we did with lower hyperspace topologies, we will denote the topology $T\wedge \tau(\mathcal{S})$ as $T(\mathcal{S})$.
\end{defn}

\begin{defn}
The collection of all subsets of $X$ that have a non-empty intersection with $S$ will be denoted as $\hit(\mathcal{S})$.
\end{defn}

\begin{prop}
Let $T$ be a miss hyperspace topology. The closure operator $\cl_{\mathcal{S}}^T$ is topological iff for each cobounded set $C$ with the property that $\downarrow C$ is a $T$-neighbourhood of $A$, we can find a cobounded set $C'\supseteq A$ such that $\downarrow C$ is a $T$-neighbourhood of $C'$.
\begin{proof}
Suppose that the closure operator $\cl_{\mathcal{S}}^T$ is topological. This means that each $\cl_{\mathcal{S}}^T(\mathcal{A})$ is closed in $T(\mathcal{S})$. Let $C$ be a cobounded set such that $\downarrow C$ is a neighbourhood of $A$. Define $S$ as the complement of $C$. Now $A$ is no element of $\cl^T(\operatorname{hit}(S)|_S)$, and therefore no element of $\cl_{\mathcal{S}}^T(\hit{S})$, and thus there is a cobounded set $C'\supseteq A$ that is not contained in $\cl^T(\operatorname{hit}(S)|_{S'})$ for some $S'\in\mathcal{S}$. The latter contains $\cl^T(\operatorname{hit}(S))$ because we are working in an upper hyperspace topology. This implies that $\downarrow C$, which is equal to the complement of $\operatorname{hit}(S)$, is a neighbourhood of $C'$.

Conversely, assume $A\not\in \cl^T(\mathcal{B}|_S)$ for some $S\in\mathcal{S}$. By definition we can find a neighbourhood $\downarrow V$ of $A$ that is disjoint from $\mathcal{B}|_{S}$. This yields that $S\setminus V$ is non-empty. Define $C$ as its complement. If for some $B\in\mathcal{B}$ we would have $B\cap S\subseteq C$, then $B\cap S$ would be a subset of $V$. Hence we obtain that $\downarrow C$ is a neighbourhood of $A$ that is disjoint from $\mathcal{B}|_S$. By assumption we can now find a $C'\supseteq A$ such that $\downarrow C$ is a neighbourhood of $C'$. This yields that $C'$ is no element of $\cl^T(\mathcal{B}|_S)$ and that $\cl_{\mathcal{S}}^T(\mathcal{B})$ is closed in $T(\mathcal{S})$, since $[\phi, X\setminus C]$ is a $\tau(\mathcal{S})$-neighbourhood of $A$ that is disjoint from $\cl_{\mathcal{S}}^T(\mathcal{B})$.
\end{proof}
\end{prop}

\begin{prop}
Let $T$ be a miss hyperspace topology. The closure operator $\cl_{\mathcal{S}}^T$ is topological iff for each $S\in\mathcal{S}$ we can find a family $(S_i)_{i\in I}$ such that $$\cl^T(\operatorname{hit}(S))=\bigcap_{i\in I}\operatorname{hit}(S_i).$$
\begin{proof}
Let $\cl_{\mathcal{S}}^T$ be a topological closure operator and take $S\in\mathcal{S}$. If $A$ is no element of $\cl^T(\operatorname{hit}(S))$, then $\downarrow C$, where $C$ is the complement of $S$, is a neighbourhood of $A$. This yields that we can find a cobounded $C'$ such that $A\subseteq C'$ and $\downarrow C$ is a neighbourhood of $C'$. This implies that, if we define $S'$ as the complement of $C'$, $\cl^T(\operatorname{hit}(S))$ is a subset of $\operatorname{hit}(S')$. The latter set clearly does not contain $A$. Since $A$ was arbitrary we obtain that the stated condition is indeed necessary.

Conversely, whenever $\downarrow C$ is a neighbourhood of a set $A$, the set $\cl^T(\operatorname{hit}(S))$, with $S$ the complement of $C$, does not contain $A$. Because $\cl^T(\hit(\mathcal{S}))$ is equal to $\bigcap _{i\in I} \hit(S_i)$ for some family $(S_i)_{i\in I}$ of bounded sets we obtain that $A$ is disjoint from a certain $S_j$. Now $S$ hits $S_j$ and thus $S_j\in \cl^T(\hit(S))$. This yields that, if we define $C'$ as the complement of $S_j$, $A\subseteq C'$ and $\downarrow C$ is a neighbourhood of $C'$ since $C'$ cannot be an element of $\cl^T(\hit(S))$.
\end{proof}
\end{prop}

\begin{prop}
If $T$ is an upper hyperspace topology, then $T(\mathcal{S})$ has the property that a set $A$ is contained in the closure of $\mathcal{A}$ iff each cobounded set that contains $A$ is in the closure of $\mathcal{A}$. We will call this property \emph{coboundedly generated}.
\begin{proof}
Let $A$ be an element of the closure of $\mathcal{A}$ and $C$ a cobounded set that contains $A$. Since we are working in an upper hyperspace topology we have that each open neighbourhood of $C$ also contains $A$. This yields that $C$ is an element of the closure of $\mathcal{A}$.

On the other hand, if $A$ is no element of the closure of $\mathcal{A}$, then we can find an open neighbourhood $\mathcal{G}$ of $A$ in $T(\mathcal{S})$ that is disjoint from $\mathcal{A}$. This open set $\mathcal{G}$ now contains a set $[S,S']$ that contains $A$ and thus $X\setminus S'$ is a cobounded set that is not contained in the closure of $\mathcal{A}$.
\end{proof}
\end{prop}

\begin{prop}
$T(\mathcal{S})$ is the finest coboundedly generated topology coarser than $T$.
\begin{proof}
Let $T^*$ be a coboundedly generated topology coarser than $T$ and let $\mathcal{G}$ be an open set in this topology. By assumption we can find a cobounded set $C$ for each element $A\in\mathcal{G}$ such that $A\subseteq C$. Because $T^*$ is an upper hyperspace topology this yields that $[\phi, X\setminus C]$ is contained in $\mathcal{G}$ and that $\mathcal{G}$ is open in $T(\mathcal{S})$.
\end{proof}
\end{prop}

\section{Bornological Convergence}

So far we have characterized the open sets in the topological reflections of both upper and lower bornological convergence structures and studied a generalization of both constructions. In what follows we wille give a description of the open sets in the topological reflection of a bornological convergence structure using uniform neighbourhoods. We will see that this complicates a possible generalization of the concept of bornological convergence to arbitrary hyperspace topologies.

\begin{defn}
A subset $G$ of a metric space is a \emph{uniform neighbourhood} of a set $A$ iff there is a $\epsilon>0$ such that $A^{\epsilon}\subseteq G$.
\end{defn}

\begin{prop}
A set $\mathcal{G}$ is open in $H(\mathcal{S})$ iff for each $A\in\mathcal{G}$ we can find $S',S\in\mathcal{S}$ and $\epsilon>0$ such that $A\in[S',S^{\epsilon}]$ and $\mathcal{G}$ is a uniform neighbourhood of $[S',S^{\epsilon}]$.
\begin{proof}
Let $A$ be an element of an open set $\mathcal{G}$ in $H(\mathcal{S})$. By definition we can find an $\epsilon>0$ and $S^*\in\mathcal{S}$ such that $\{B|A\cap S^*\subseteq B^{\epsilon} \textrm{ and } B\cap S^*\subseteq A^{\epsilon}\}$. Define $S'$ and $S$ respectively as $A\cap S^*$ and $S^*\setminus A^{\epsilon}$. Clearly we have that $A$ is contained in $[S',S^{\epsilon}]$. Take $C\in[S',S^{\epsilon}]$ and a set $B$ that has a Hausdorff distance less than $\epsilon$ from $C$. First of all, we have $A\cap S^*\subseteq C$ and $C\subseteq B^{\epsilon}$. Moreover, we have $B\cap S^*\subseteq C^{\epsilon}\cap S^*$ and $C^{\epsilon}\cap S^*\subseteq A^{\epsilon}$. All this yields that $B$ is an element of $\mathcal{G}$ and that $\mathcal{G}$ is a uniform neighbourhood of $[S',S^{\epsilon}]$.

Conversely, assume that $A$ is contained in a set $[S',S^{\epsilon}]$ of which $\mathcal{G}$ is a uniform neighbourhood. By definition we can find an $\epsilon'>0$ such that all sets that have a Hausdorff distance less than $\epsilon'$ from an element of $[S',S^{\epsilon}]$ are contained in $\mathcal{G}$. Let $\delta>0$ be the minimum of $\epsilon$ and $\epsilon'$ and define $S^*$ as $S\cup S'$. Suppose that we have a set $B$ that satisfies $A\cap S^*\subseteq B^{\delta}$ and $B\cap S^*\subseteq A^{\delta}$. This yields that $B\cup S'$ is an element of $[S',S^{\epsilon}]$ and that $H(B,B\cup S')\leq \epsilon'$. Hence we obtain that $B$ is an element of $\mathcal{G}$ and that the latter is an open set in $H(\mathcal{S})$.
\end{proof}
\end{prop}

The previous proposition shows that the uniform structure generated by the Hausdorff distance plays a crucial role in the construction of $H(\mathcal{S})$. This suggests that this particular concept doesn't translate as easily to general hyperspace topologies as was the case with lower and upper bornological convergence structures. In both cases we found that the topological reflection of such a bornological modification was equal to the infimum of the original topology $T$ and the topology $\tau(\mathcal{S})$. The modification of $H$ that leads to $H(\mathcal{S})$, the topological reflection of $\cl_{\mathcal{S}}^H$, doesn't seem to have this advantage. The next example shows that $H(\mathcal{S})$ is in general not equal to $H\wedge \tau(\mathcal{S})$.

\begin{exam}
Let $X$ be the metric space $[-1,1]$ endowed with the Euclidian metric and $\mathcal{S}$ the ideal of finite subsets of $X$. Define $\mathcal{B}$ as the collection of all sets that are either of the form $\{-\delta, \delta\}$ with $\delta\in [0,1]$ or are equal to the empty set. It is clear that $\mathcal{B}$ is closed for the Hausdorff distance on the hyperspace of $X$. Now let $A$ be an element in the $\tau(\mathcal{S})$-closure of $\mathcal{B}$. If $A$ is the empty set, then it is contained in $\mathcal{B}$ by definition. Suppose that $A$ is non-empty and take $\delta\in A$. Since $\{-\delta, \delta\}$ is an element of $\mathcal{S}$ we should be able to find a $B\in \mathcal{B}$ such that $A\cap \{-\delta, \delta\}$ equals $B\cap \{-\delta, \delta\}$. This can only be true if $B$ is equal to $\{-\delta, \delta\}$ and $A$ contains $\{-\delta, \delta\}$. Assume that $\epsilon$ is an element of $A$ that is not an element of $\{-\delta, \delta\}$ and define $S$ as $\{-\delta, \epsilon, \delta\}$. Again we can find a $B\in \mathcal{B}$ such that $A\cap S=B\cap S$. It is clear that this is not possible. Hence $A$ is equal to $\{-\delta, \delta\}$ and therefore contained in $\mathcal{B}$. This means that $\mathcal{B}$ is closed in $\tau(\mathcal{S})$.

To prove that $\mathcal{B}$ is not closed in $H(\mathcal{S})$ we will show that $\{1\}$ is an element of the bornological closure of $\mathcal{B}$. Let $S$ be a finite subset of $[-1,1]$ and take $\epsilon>0$. Choose a $\delta\in ]1-\epsilon, 1]$ such that $\delta\not\in S$ and denote $\{-\delta, \delta\}$ as $B$. We then obtain that $\{1\}\cap S\subseteq B^{\epsilon}$ and $B\cap S\subseteq \{1\}^{\epsilon}$. This means that $\{1\}$ is an element of the bornological closure of $\mathcal{B}$. It is, however, by definition no element of $B$ itself.
\end{exam}

Like we did before with the upper and lower bornological convergences we would like to know if the closure operator of the topology $H(\mathcal{S})$ is equal to $\cl_{\mathcal{S}^*}^H$ for some ideal $\mathcal{S}^*$. Example \ref{H(S)bornological} shows that this need not always be the case.

\begin{prop}
If the closure operator associated with the topology $H(\mathcal{S})$ coincides with $\cl_{\mathcal{S}^*}^H$ for some ideal $\mathcal{S}^*$, then $\mathcal{S}^*$ equals $\mathcal{S}^+$.
\begin{proof}
By assumption we get that $\cl_{\mathcal{S}}^H$ is finer than $\cl_{\mathcal{S}^*}^H$. In \cite{Lechicki} we find that this implies that $\mathcal{S}^*\subseteq \mathcal{S}$ and that $\mathcal{S}^*$ is closed under small enlargements. Hence we obtain that $\mathcal{S}^*$ is contained in $\mathcal{S}^+$. On the other hand we have that $\cl_{\mathcal{S}^+}^H$ is topological and thus coarser than $\cl_{\mathcal{S}^*}^H$. This yields $\mathcal{S}^+\subseteq \mathcal{S}^*$ and thus both ideals are equal.
\end{proof}
\end{prop}

\begin{exam}\label{H(S)bornological}
Let $\mathbb{R}^2$ be the real plane with the Euclidian metric and $\mathcal{S}$ the bornology of finite subsets of $\mathbb{R}^2$. Define $P$ as $\{(x,y)\in\mathbb{R}^2| y>0\}$. The set $\hit(P)$ is open in $H(\mathcal{S})$. Take $A\in \hit(P)$ and $(x,y)\in A$ with $y>0$. If we define $\epsilon$ as $y/2$ and $S$ as $\{(x,y)\}$, then $\{B|A\cap S\subseteq B^{\epsilon} \textrm{ en } B\cap S\subseteq A^{\epsilon}\}$ is a neighbourhood of $A$ that is subset of $\hit(P)$. The ideal $\mathcal{S}^+$, however, only contains the empty set. This means that $\cl_{\mathcal{S}^+}^H$ is equal to the closure operator of the indiscrete topology. 
\end{exam}

\end{document}